\documentclass{article}
\usepackage{graphicx} % Required for inserting images
\usepackage{amsfonts,amssymb,amsmath,amsthm}
\usepackage{xcolor}
\usepackage{appendix}
\usepackage{nicefrac}
\usepackage{placeins}
\usepackage{tikz,pgfplots}
\usetikzlibrary{shapes.geometric, arrows}
\numberwithin{equation}{section}

\allowdisplaybreaks

\newcommand{\pdet}{{\det}^\dagger}
\newcommand{\D}{D}

\tikzstyle{box} = [rectangle, rounded corners, minimum width=4cm, minimum height=.7cm,text centered, draw=black, fill=black!30]
\tikzstyle{midbox} = [rectangle, rounded corners, minimum width=2.62cm, minimum height=.7cm,text centered, draw=black, fill=black!30]
\tikzstyle{arrow} = [thick,->,>=stealth]

\title{Reprojection methods for Koopman-based modelling and prediction\thanks{K.~Worthmann gratefully acknowledges funding by the Deutsche Forschungsgemeinschaft (DFG, German Research Foundation) -- {Project-ID 507037103}}}
\author{Pieter van Goor$^1$, Robert Mahony$^1$, Manuel Schaller$^2$, and Karl Worthmann$^2$}
\date{%
    \normalsize
    $^1$Systems Theory and Robotics Group, School of Engineering, Australian National University, Canberra, Australia\\
    $^2$Optimization-based Control Group, Institute of Mathematics, Technische Universit\"at Ilmenau, Germany\\[2ex]%
    July 2023
}

\newtheorem{theorem}{Theorem}

\newtheorem{proposition}[theorem]{Proposition}

\newtheorem{definition}[theorem]{Definition}
\newtheorem{remark}[theorem]{Remark}

\newtheorem{example}[theorem]{Example}

\begin{document}

\maketitle

\begin{abstract}
Extended Dynamic Mode Decomposition (eDMD) is a powerful tool to generate data-driven surrogate models for the prediction and control of nonlinear dynamical systems in the Koopman framework. 
	In eDMD a compression of the lifted system dynamics on the space spanned by finitely many observables is computed, in which the original space is embedded as a low-dimensional manifold.
	While this manifold is invariant for the infinite-dimensional Koopman operator, this invariance is typically not preserved for its eDMD-based approximation.
	Hence, an additional (re-)projection step is often tacitly incorporated to improve the prediction capability. 
	We propose a novel framework for consistent reprojectors respecting the underlying manifold structure. 
	Further, we present a new geometric reprojector based on maximum-likelihood arguments, which significantly enhances the approximation accuracy and preserves known finite-data error bounds. 
\end{abstract}

\section{Introduction}\label{sec:intro}

In the Koopman framework nonlinear dynamical systems are lifted to the infinite-dimensional space of observables, in which the system dynamics are governed by a semi-group of linear operators. 
Since a compression of the Koopman operator can be efficiently computed in a purely data-based manner based on the extended Dynamic Mode Decomposition (eDMD), see, e.g., \cite{WillMatt15}, Koopman-based prediction and control has attracted considerable attention in recent years, see, e.g., the collection~\cite{MaurMezi20}, the recent survey~\cite{BevaSosnHirc21}, and the references therein.
In eDMD, finitely-many observables are evaluated along a finite number of sample trajectories to compute a compression of the infinite-dimensional Koopman operator by means of a regression problem~\cite{KlusNusk20}. 
The approximation is subject to an estimation error due to a finite amount of data~\cite{NuskPeit23}, and a projection error stemming from a finite dictionary size~\cite{SchaWort22}.

\noindent The observables map the state space to a manifold in the lifted space, which is preserved by the Koopman operator, i.e., both the lifted state trajectories and the corresponding Koopman flow evolve on this manifold. 
However, the finite-dimensional approximation does not satisfy this property if one neglects particular cases typically linked to Koopman invariance of the dictionary~\cite{BrunBrun16}, see also the recent work~\cite{Mezi20}. While the mentioned references provide, at least up to a certain degree, a remedy for the prediction and analysis of dynamical systems, the respective conditions w.r.t.\ control systems are quite restrictive except for the drift-free case~\cite{GoswPale21}. Hence, we are concerned with the case that the dictionary is not Koopman invariant, which is often present in practice. This is of paramount importance since the %underlying reasoning is that 
learned compression is based on information on the manifold only and, thus, may exhibit large errors if applied on the lifted space, but not on the manifold itself. The respective errors often even amplify for increasing dictionary size increases. %--~in contrast to the projection error analysis.
Often a (re-) projection step is tacitly incorporated to ensure consistency, i.e., projecting back to the manifold~\cite{MaurGonc16}, and, thus, to counteract these deteriorating effects. To be more precise, following each prediction step, the propagated observables are projected back onto the manifold before the surrogate model is used for the next prediction step. If the coordinate functions are included in the dictionary, the canonical choice is the coordinate projection~\cite{TeraShir21}. 
Whereas the coordinate projection may be an attractive choice due to its simplicity, we demonstrate that it is, in general, by far not the best. 
Additionally, it requires that the coordinate functions are included in the dictionary.

%Karl: contribution statement
%In this work, we present a framework for closest-point projections in a particular metric on the ambient space of the manifold. %We show, that the coordinate projection is a closest-point projection in a particular metric.
The contribution of this paper is twofold. 
First, we %rigorously 
prove that the additional projection step cannot deteriorate the overall performance %too 
much. To this end, we link the respective error to the one resulting from the regression problem for computing the finite-dimensional compression of the Koopman operator. In addition, we explicate why the coordinate projection performs well %essentially does the job 
if the right-hand side of the differential equation is contained in the %span of the observables $\matbb{V}$.
dictionary~$\mathbb{V}$.
Second, we propose a novel framework for closest-point projections containing the \textit{standard} coordinate projection as well as an alternative \textit{geometric} projection based on a Maximum-Likelihood estimator. 
To this end, we introduce a semi-inner product on the ambient space $\mathbb{R}^N$, where $N$ is the dimension of the dictionary $\mathbb{V}$, which induces a Riemannian metric on the manifold.
This allows for different weightings, and we show that this does not interfere with solving the %previously-solved 
regression problem. 
To be more precise, we prove that the solution of the $L^2$-regression problem is also a solution w.r.t.\ this new weighted counterpart. %$L^2$-problem. 
In conclusion, the different projections correspond to different choices of semi-inner product. 
%Furthermore, we present numerical simulations showing that this novel approach is superior w.r.t.\ its approximation accuracy by considering the one-step prediction error.
To illustrate the superiority of the presented approach w.r.t.\ approximation accuracy, we provide several numerical examples.

%Karl: Outline
The outline of the paper is as follows: In Section~\ref{sec:edmd}, we briefly recap eDMD in the Koopman framework and in Section~\ref{sec:problem_formulation} we provide the problem formulation. Then, in Section~\ref{sec:coordinate_projection}, we consider the coordinate projection to show that a projection step between predictions in the lifted space is (highly) %(very) 
beneficial. In Section~\ref{sec:closest_point_projection}, the novel projection framework is introduced and the key results are presented before the geometric projection and numerical simulations are conducted in the following two sections. In Section~\ref{sec:conclusions}, conclusions are drawn before a brief outlook is given.

\noindent \textbf{Notation}: Let $\mathbb{R}$ be the field of the real numbers. Further, for integers $a,b \in \mathbb{Z}$ with $a \leq b$, we set $[a:b] := [a,b] \cap \mathbb{Z}$. For $\eta > 0$, $\mathcal{B}_\eta(\hat{x})$ denotes the ball centered at~$\hat{x} \in \mathbb{R}^d$ with radius~$\eta$, i.e., the set $\{ x \in \mathbb{R}^d: \| x - \hat{x} \| < \eta \}$, where $\|\cdot\| : \mathbb{R}^d \rightarrow \mathbb{R}$ is the Euclidean norm. Moreover, for two sets~$A, B \subset \mathbb{R}^d$, $A \oplus B := \{ x + y : x \in A, y \in B \}$ is the Pontryagin sum. 
The pseudodeterminant $\pdet(W)$ of a matrix $W \in \mathbb{R}^{N\times N}$ is defined as the product of the nonzero singular values of~$W$.
For a function $\varphi: \mathbb{R}^n \to \mathbb{R}^m$, the differential at a point $x \in \mathbb{R}^n$ is written as $\mathrm{D} \varphi(x) \in \mathbb{R}^{m\times n}$.

\section{Extended Dynamic Mode Decomposition} 
\label{sec:edmd}

For a compact set~$\mathbb{X}$ and (sufficiently large) $\eta > 0$, we consider the system dynamics
\begin{equation}\label{eq:ode}
\dot{x}(t) = f(x(t))
\end{equation}
with Lipschitz continuous vector field $f: \D := \mathbb{X} \oplus \mathcal{B}_\eta(0) \subset \mathbb{R}^d \rightarrow \mathbb{R}^d$. We denote the solution of~\eqref{eq:ode} at time $t \in \mathbb{R}_{\geq 0}$ for the initial condition $x(0) = \hat{x} \in \mathbb{X}$ by $x(\,\cdot\,;\hat{x})$ on its maximal interval $[0,t_{\hat{x}})$ of existence. For given $\Delta t > 0$, $t_{\hat{x}} \geq \Delta t$ holds for all $\hat{x} \in \mathbb{X}$ if $\eta > 0$ is sufficiently large, which is tacitly assumed in the following to streamline the presentation. 
Alternatively, forward invariance of the set~$\mathbb{X}$ w.r.t.\ the flow of the dynamical system governed by~\eqref{eq:ode} is assumed, see, e.g., \cite{ZhanZuaz23}.

\begin{figure}[htb]
	\centering
	\resizebox{.45\columnwidth}{!}{%	
		\begin{tikzpicture}[node distance=1.4cm,shorten <=.25cm, shorten >=.25cm]
		\node (start) [midbox,text width=2.5cm] {observable \\$\varphi\in L^2(D,\mathbb{R}^n)$}; % \in {\X}$};
		\node (startleft)[left of=start, node distance = .4cm] {};
		
		\node (startright)[right of=start, node distance = .5cm] {};	
		\node (startrightarrow)[below of=startright, node distance = .15cm] {};
		\node (mid) [midbox, below of=start] {$\mathcal{K}^t \varphi$};
		\node (midleft)[left of=mid, node distance = .4cm] {};
		\node (midright)[right of=mid, node distance = .5cm] {};
		%\draw [arrow] (startleft) -- node[] {} (midleft);
		\draw [arrow] (startrightarrow) --  node[] {} (midright);
		\node (dummy) [below of = startleft, node distance = .75cm] {\small Koopman};
		\node (bot) [midbox, below of=mid] {$(\mathcal{K}^t \varphi)(\hat{x})$};
		\node (botleft)[left of=bot, node distance = .4cm] {};
		\node (botright)[right of=bot, node distance = .5cm] {};
		%\draw [arrow] (midleft) -- node[] {} (botleft);
		\draw [arrow] (midright) -- node[] {} (botright);	    
		\node (dummy) [below of = midleft, node distance = .65cm] {\small evaluate};
		\node (start2) [midbox, left of=start,node distance=3.3cm,text width=2.5cm] {initial state\\ $\hat{x}\in \mathbb{X}$}; % \in L^2({\X})$};
		\node (start2left)[left of=start2, node distance = .5cm] {};
		\node (start2leftarrow)[below of=start2left, node distance = .12cm] {};    
		\node (start2right)[right of=start2, node distance = .4cm] {};
		\node (mid2) [midbox, below of=start2] {$x(t;\hat{x})$};
		\node (mid2left)[left of=mid2, node distance = .5cm] {};
		\node (mid2right)[right of=mid2, node distance = .4cm] {};
		\draw [arrow] (start2leftarrow) -- node[] {} (mid2left);
		%\draw [arrow] (start2right) -- node[] {} (mid2right);
		\node (dummy2) [below of = start2right, node distance = .75cm] {\small ODE};
		\node (bot2) [midbox, below of=mid2] {$\varphi(x(t;\hat{x}))$};
		\node (bot2left)[left of=bot2, node distance = .5cm] {};
		\node (bot2right)[right of=bot2, node distance = .4cm] {};
		\draw [arrow] (mid2left) -- node[] {} (bot2left);
		%\draw [arrow] (mid2right) -- node[] {} (bot2right);	    
		\node (dummy) [below of = mid2right, node distance = .65cm] {\small evaluate};
		\node (equation) [right of=bot2,node distance = 1.65cm] {};
		\node (end)[midbox,below of=equation]{$\varphi(x(t;\hat{x}))=(\mathcal{K}^t\psi)(\hat{x})$};
		\node (asd)[below of=equation,node distance =.65cm]{equate};
		\node(startarrowleft) [below of=botleft, node distance=.2cm] {};
		\node(startarrowright) [below of=bot2right, node distance=.2cm] {}; 
		\draw [arrow] (startarrowleft) -- node[] {} (end);
		\draw [arrow] (startarrowright) -- node[] {} (end);
		\end{tikzpicture}}
	\caption{Schematic sketch of the Koopman framework.}% The evaluation is conducted w.r.t.\ time and space resp.}
	\label{fig:sketch}
\end{figure}
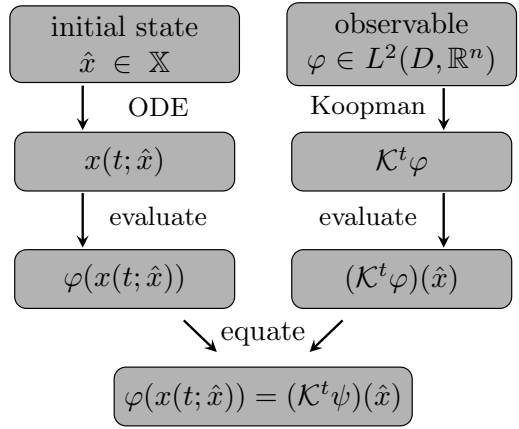

The Koopman %operator 
semigroup~$(\mathcal{K}^t)_{t\geq 0}$ %is a semigroup 
of linear operators on $L^2(\D,\mathbb{R})$ is
defined via the identity
\begin{equation}\label{eq:Koopman}
(\mathcal{K}^t \varphi)(\hat{x}) = \varphi(x(t;\hat{x})) \qquad\forall\,\hat{x} \in \mathbb{X}, \varphi \in L^2(\D,\mathbb{R})
\end{equation}
for all $t < t_{\hat{x}}$ and, thus, in particular on the interval $[0,\Delta t]$.
By means of this semigroup, one may either propagate the observable~$\varphi$ forward in time using the Koopman operator~$\mathcal{K}^t$ and evaluate the propagated observable at~$\hat{x}$ or evaluate the observable~$\varphi$ at the solution $x(t;\hat{x})$ as depicted in Figure~\ref{fig:sketch}.

\noindent For $N \in \mathbb{N}$ and linearly independent observables $\psi_i \in L^2(\D,\mathbb{R})$, $i \in [1:N]$, $\mathbb{V} := \operatorname{span} \{ \psi_i \mid i \in [1:N] \}$ is called the dictionary. 
Define the vector-valued function
\begin{align}\label{eq:Psi}
	\Psi:\D \to \mathbb{R}^N, \qquad x \mapsto (\psi_1(x)\,\ldots\,\psi_N(x))^\top.
\end{align}
Invoking the identity~\eqref{eq:Koopman}, we get
$$
\mathcal{K}^t \Psi = \begin{pmatrix}
\mathcal{K}^t \psi_1\\
\vdots\\
\mathcal{K}^t\psi_N
\end{pmatrix} = \Psi(x(t;\hat{x})) \in L^2(\mathbb{X},\mathbb{R})^N \simeq L^2(\mathbb{X},\mathbb{R}^N).
$$
Then, the approximation of the Koopman operator is defined as a best-fit solution in an $L^2$-sense by means of the regression problem
\begin{align}\label{eq:fitting}
	\begin{split}
		\hat{K}
		& = \operatorname{arg} \hspace*{-0.5cm}\min_{K\in \mathbb{R}^{N \times N} \hspace*{0.4cm}}\hspace*{-0.2cm} \int_\mathbb{X} \|\Psi(x(t;\hat{x})) - K\Psi(\hat{x})\|^2_2 \,\mathrm{d}\hat{x}.
	\end{split}
\end{align}
A solution to Problem~\eqref{eq:fitting} can easily be calculated by 
\begin{align}\label{eq:KI}
	\hat{K} = \Psi_Y \Psi_X^{-1},
\end{align}
where $\Psi_X$ and $\Psi_Y$ are given by $\Psi_X = \int_\mathbb{X} \Psi(\hat{x})\Psi(\hat{x})^\top\,\mathrm{d}\hat{x}$ and $\Psi_Y = \int_\mathbb{X}\Psi(x(t;\hat{x}))\Psi(\hat{x})^\top \,\mathrm{d}\hat{x}$, respectively. $\hat{K}$ is a compression of $\mathcal{K}^t$, that is, $\hat{K} = {P}_\mathbb{V}\mathcal{K}^t_{\vert\mathbb{V}}$ holds,
where $P_\mathbb{V}$ is the $L^2$-orthogonal projection onto $\mathbb{V}$. The %respective 
projection error was first analyzed in~\cite{ZhanZuaz23} by means of a dictionary of finite elements, see also~\cite{SchaWort22} for an extension to control systems.
In particular, if $\mathbb{V}$ is Koopman-invariant, then $\hat{K} = \mathcal{K}^t_{\vert \mathbb{V}}$.
\begin{remark}\label{rem:finite_data}
	An empirical estimator of $\hat{K}$ given $m\in \mathbb{N}$ i.i.d.\ data points $x_1,\ldots,x_m$ can be computed via
	$$
	\hat{K}_{m} = \operatorname{arg} \hspace*{-0.5cm}\min_{K\in \mathbb{R}^{N \times N} \hspace*{0.4cm}}\hspace*{-0.2cm} \sum_{i=1}^m \|\Psi(x(t;x_i)) - K \Psi(x_i)\|_2^2.
	$$
	For $m \geq N$, $\hat{K}_m = (\Psi_X^m {\Psi_Y^m}^\top)(\Psi_X^m {\Psi_X^m}^\top)^{-1}$ is a closed-form solution
	using the $(N \times m)$-data matrices $\Psi_X^m$ and $\Psi_Y^m$ with entries $\psi_i(x_j)$ and $\psi_i(x(t;x_j))$, $(i,j) \in [1:N] \times [1:m]$, respectively. 

	The convergence $\hat{K}_m \rightarrow \hat{K}$ for $m \rightarrow \infty$ follows by the law of large numbers~\cite[Section 4]{KordMezi18}. For finite-data error bounds we refer to~\cite{NuskPeit23} and the references therein, where also an extension to Stochastic Differential Equations (SDEs) with ergodic sampling (along a single, sufficiently long trajectory) is given. For a recent result in reproducing kernel Hilbert spaces, we refer to~\cite{PhilScha23}.
\end{remark}

In fact, the regression problem~\eqref{eq:fitting} is the ideal formulation for the compression in the infinite-data limit. In our numerical simulations, we solve an approximation of this using $m = 10,000$ i.i.d.\ data points drawn from the compact set~$\mathbb{X}$. %, cp.~Remark~\ref{rem:finite_data} for details.

\section{Problem formulation}\label{sec:problem_formulation}
The dictionary~$\mathbb{V}$ is defined by the span of linearly-independent observables $\psi_1,...,\psi_N$.
In \eqref{eq:fitting}, $\hat{K}$ is computed by stacking the observables $\psi_1,..,\psi_N$ into $\Psi \in L^2(\mathbb{X}, \mathbb{R})^N$. Correspondingly, we define the set 
\begin{align}\label{eq:M}
	M := \operatorname{im}(\Psi) = \{ \Psi(x) \mid x \in \mathbb{X} \} \subset \mathbb{R}^N.
\end{align}
By definition, the set~
$M$ is invariant w.r.t.\ the Koopman operator~$\mathcal{K}_t$, $t \in [0,\Delta t]$, i.e., 
%By definition, the Koopman operator preserves the manifold~$M$, i.e., 
\begin{equation}\label{eq:Koopman_manifold}
(\mathcal{K}^t \Psi)(\hat{x}) = \Psi(x(t;\hat{x})) \in M \qquad\forall\,\hat{x} \in \mathbb{X}.    
\end{equation}
The following example taken from \cite{BrunBrun16} nicely illustrates that this property also holds for the eDMD-based surrogate model if the dictionary~$\mathbb{V}$ is Koopman-invariant, i.e., $\mathcal{K}^t\mathbb{V} \subseteq \mathbb{V}$. 
\begin{example}\label{ex:invariance}
	Consider the system
	\begin{align}\label{eq:rhs_ex1}
		\frac{\mathrm{d}}{\mathrm{d}t}\begin{pmatrix}
			x_1(t)\\x_2(t)
		\end{pmatrix}
		=
		\begin{pmatrix}
			x_1(t)\\ \lambda(x_2(t) - x_1^2(t))
		\end{pmatrix}.
	\end{align}
	in $\mathbb{R}^2$ with $\lambda \in \mathbb{R}$.
	Choosing the observables $\psi_1(x) = x_1$, $\psi_2(x) = x_2$, and $\psi_3(x) = x_1^2$, we get $\dot{\psi}_1(x(t)) = \psi_1(x(t))$, $\dot{\psi}_2(x(t)) = \lambda (\psi_2(x(t)) - \psi_3(x(t)))$, $\dot{\psi}_3(x(t)) = 2 \psi_3(x(t))$, which can be written as the linear system 
	\begin{align}\label{eq:linear_example}
		\frac{\mathrm{d}}{\mathrm{d}t} \begin{pmatrix}
			y_1\\y_2\\y_3 
		\end{pmatrix}(t)
		= 
		\begin{pmatrix}
			1&0&0\\
			0&\lambda&-\lambda\\
			0&0&2
		\end{pmatrix}
		\begin{pmatrix}
			y_1\\y_2\\y_3 
		\end{pmatrix}(t) =: Ay(t).
	\end{align}
	As the prediction of observables in $\mathbb{V} = \operatorname{span}\{ \psi_1,\psi_2,\psi_3\}$ can be performed by means of \eqref{eq:linear_example}, 
	the subspace~$\mathbb{V}$ is Koopman invariant and the matrix representation $K^t \in \mathbb{R}^{3\times 3}$ of the Koopman operator $\mathcal{K}^t|_{\mathbb{V}}$ w.r.t.\ the basis $\{\psi_1,\psi_2,\psi_3\}$ is given by $K^t = e^{tA}$.
	The prediction of an observable $\varphi = \sum_{i=1}^3 a_i \psi_i$ along the flow emanating from~$\hat{x}$ is given by $\varphi(x(t;\hat{x})) =
	\langle a, e^{tA} \Psi(\hat{x}) \rangle_2$. 
\end{example}

The invariance of~$\mathbb{V}$ for Example~\ref{ex:invariance} is preserved if the $x_1$-component of the right-hand side is multiplied by $\mu \in \mathbb{R} \setminus \{0\}$ or the term~$x_1^2$ in the $x_2$-component is replaced by an arbitrary polynomial $p(x_1)$, see~\cite{BrunBrun16}.
However, this desirable property does not hold in general as showcased in the following example.
\begin{example}
	Let us replace the linear term~$x_1$ in the first component of Example~\ref{ex:invariance} by~$-x_1^2$, i.e.,  
	\begin{align*}
		\frac{\mathrm{d}}{\mathrm{d}t}\begin{pmatrix}
			x_1(t)\\x_2(t)
		\end{pmatrix} = \begin{pmatrix}
			-x_1^2(t)\\ \lambda(x_2(t) - x_1^2(t))
		\end{pmatrix}.
	\end{align*}
	Then, the dictionary spanned by $\psi_i(x) = x_i$, $i \in \{1,2\}$, and $\psi_i(x) = x_1^{i-1}$, $i \in \mathbb{N}_{\geq 3}$, is Koopman invariant, but infinite dimensional.
\end{example}

\begin{figure}[htb]
	\centering
	\begin{tikzpicture}[scale = 1.25]
	% Help lines
	%  \draw[help lines] (-1,-1) grid (8,6);
	
	% Manifold
	\draw[smooth, tension=0.5] plot coordinates{(1,2) (2,2.4) (4,1.78) (6.3,2.8)}{};
	\draw[smooth, tension=0.5] plot coordinates{(6.3,3.8) (4.5,3) (2,3.5) (1.1,3)}{};
	\draw[smooth, tension=1] plot coordinates{(1,2) (1.1,2.5) (1.1,3)}{};
	\draw[smooth, tension=1] plot coordinates{(6.3,2.8) (6.4,3.4) (6.3,3.8)}{};
	%\draw[smooth cycle, tension=0.5] plot coordinates{(1,2) (2,2.5) (4,2) (6,3) (6.5,4) (4.5,3) (2,3.5) (1,3)}{};
	\node at (6.5,4.1) {$\Psi(\mathbb{R}^d)$};
	%\draw[line width=.08em,smooth cycle, tension=0] plot coordinates{(0,0) (0.8,1) (6.8,1) (6,0)} node[label=right:$\mathbb{R}^d$] {};
	\draw[line width=.08em] (0.5,0)--(1,1)--(7,1)--(6.5,0) -- cycle;%node[label=right:$\mathbb{R}^d$] {};
	\node (a) at (7,0.2) {$\mathbb{R}^d$};
	
	\draw[line width=.08em,->] (0.5,0) -- (0.5, 4) node [label=right:$\mathbb{R}^{N-d}$] {};
	
	% Stuff in Rd
	\draw[smooth, tension=0.5] plot coordinates{(2,0.2) (2.5,0.5) (3.5,0.2) (4,0.4)} node[above left=-.05cm] {$x(\Delta t;\hat{x})$};
	\draw[thin, fill=black] (2,0.2) circle (1.5pt) node[label=left:$\hat{x}$]{};
	\draw[thin, fill=black] (4,0.4) circle (1.5pt);
	
	% Stuff on Manifold $(\mathcal{K}^t \psi)(\hat{x})$
	\draw[smooth, tension=0.5] plot coordinates{(2,2.8) (2.5,2.9) (3.5,2.3) (4,2.5)} node[above=-.01cm] {$\Psi(x(\Delta t;\hat{x}))$};
	\draw[thin, fill=black] (2,2.8) circle (1.5pt) node[label=above:$\Psi(\hat{x})$]{};
	\draw[thin, fill=black] (4,2.5) circle (1.5pt);
	\draw [dotted] (2,2.8) -- (2,0);
	\draw [dotted] (4,0.4) -- (4,2.5);
	\node (b) at (3,2.35) [rotate=-32] {$(\mathcal{K}^t\Psi(\hat{x}))_{t\geq 0}$};
	
	% Koopman approximation
	%\draw [line width=.01em, dashed,red] (2,2.8) -- (5,1.5); %node[label=right:$\hat{K}\psi(\hat{x})$]{};
	\draw[thin,red,fill=red] (5,1.5) circle (1.5pt);
	\node (a) at (5.55,1.5) {\textcolor{red}{$\hat{K}\Psi(\hat{x})$}};
	% Projection onto manifold
	\draw [dashed] (4.98,1.55) -- (4.7,2.6) ;%node[label=above right:$\pi_W(\hat{K}\Psi(\hat{x}))$]{};
	\draw[thin, fill=black] (4.7,2.6) circle (1.5pt);
	\node (a) at (7,0.2) {$\mathbb{R}^d$};
	\node (b) at (5.5,2.9) {$\pi_W\!\circ\!\hat{K}\Psi(\hat{x})$};
	% Simple brace
	% \draw [decorate, decoration = {brace,raise=.3em}] (4,0.4)--(4.7,0.8) node[left=.55cm]{coordinate error};
	% Projection onto Rd
	\draw [dashed] (4.7,2.6) -- (4.7,0.8);
	\node (c) at (5.5,0.5){$\Psi^{-1}\!\circ\!\pi_W\!\circ\!\hat{K}\Psi(\hat{x})$};
	\draw[thin, fill=black] (4.7,0.8) circle (1.5pt);
	\end{tikzpicture}
	\caption{Geometric projection after applying the approximation~$\hat{K}$.}
	\label{fig:manifold}
\end{figure}
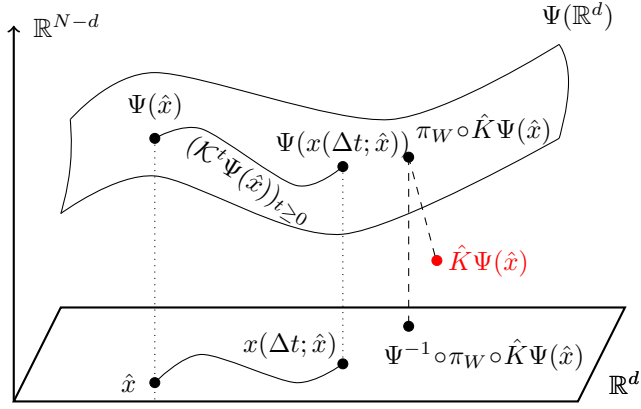

In conclusion, one cannot expect %that 
$\hat{K}\Psi(\hat{x}) \in M$ %holds 
for the approximated Koopman operator~$\hat{K}$ as depicted in Figure~\ref{fig:manifold}.
This causes two issues in using $\hat{K}$ (or data-driven approximations thereof) to model the system dynamics~\eqref{eq:ode}.
The first is that it is unclear how to recover the state values underlying the propagated observables.
Specifically, if $\hat{K}\Psi(\hat{x}) \notin M$, then by definition there is no value $x \in D \subset \mathbb{R}^d$ satisfying $\Psi (x) = \hat{K} \Psi(\hat{x})$.
The second issue is that the learning process, i.e., the regression problem~\eqref{eq:fitting}, only uses measurements of the form $z = \Psi(x)$, i.e, only points contained in the set~$M$ %on the manifold 
are taken into account.
Hence, one cannot expect $\hat{K}z$ to be meaningful if $z \notin M$, which may render a repeated application of~$\hat{K}$ questionable. 
Both of these issues can be mitigated by \emph{projecting} the dynamics $z^+ = \hat{K}z$, $z = \Psi(x)$, \emph{back} to the set~$M$ after each iteration, see, e.g.,~\cite{MaurGonc16}.

Within this paper, we propose a framework for understanding a wide class of projections and propose two particular choices: the regularly used \emph{coordinate projection} and our newly introduced \emph{geometric projection}. Such a projection step in particular is crucial for future applications in Koopman-based (predictive) control \cite{SchaWort22,BoldGrun23} using \mbox{eDMDc}~\cite{KordMezi18b} or a bilinear surrogate model~\cite{WillHema2016}, where the construction of Koopman-invariant subspaces is a highly-nontrivial issue, see, e.g., \cite{GoswPale21}. % in Sections~\ref{sec:coordinate_projection} and~\ref{sec:geometric_projection}, respectively.

\section{The coordinate projection}\label{sec:coordinate_projection}
In this section, we consider the coordinate projection --~an approach that has been used by various authors, including those developing neural-network based EDMD~\cite{TeraShir21}.
If the first $d$ observables are chosen to be the coordinate functions, i.e., $\psi_i(x) = x_i$ holds for all $i \in [1:d]$, we get
\begin{align}\nonumber %\label{eq:Psi}
	\Psi(x) = \begin{pmatrix}
		x \\ \bar{\Psi}(x)
	\end{pmatrix} := \big( x_1\,\ldots\,x_d\ \psi_{d+1}(x)\,\ldots\,\psi_N(x) \big)^\top,
\end{align}
where $\bar{\Psi} = \Psi_{[d+1:N]}$ consists of the last $N-d$ components of $\Psi$. Then, assuming that $\bar{\Psi}: \mathbb{R}^d \rightarrow \mathbb{R}^{N-d}$ is a smooth function, the set~$M$ is a graph and, thus, a smooth manifold. %\RM{[The smooth manifold was already claimed. It is the graph property that matters]}
We emphasize that, for each $z \in M$, there exists a unique $x = z_{[1:d]} \in \mathbb{R}^d$ such that $\Psi(x) = z$ holds, i.e., $\Psi$ is invertible by simply taking the first $d$ coordinates.
Hence, the coordinate projection $\pi : \mathbb{R}^N \to M$ is defined by 
\begin{align}\label{eq:coordinate_projection}
	\pi(z) := (z_{1,...,d}, \bar{\Psi}(z_{1,...,d}))\in M,
\end{align}
for all $z \in \mathbb{R}^N$.
The associated approximated discrete-time dynamics on $\mathbb{X} \subset \mathbb{R}^d$ are defined as
\begin{align}\label{eq:koopman_approx_dynamics}
	x^+ = \hat{F}_{\pi}(x) := \Psi^{-1} \circ \pi(\hat{K} \Psi(x)).
\end{align}
Using the particular structure of the coordinate projection~\eqref{eq:coordinate_projection}, this simplifies to
\begin{align*}
	\hat{F}_{\pi}(x) % &= \Psi^{-1} \circ \pi_D(\hat{K} \Psi(x)), \\
	&= \Psi^{-1} \circ \pi \begin{pmatrix}
		\hat{K}_{[1:d]} \Psi(x)\\ \hat{K}_{[d+1:N]} \Psi(x)
	\end{pmatrix} 
	= \Psi^{-1} \begin{pmatrix}
		\hat{K}_{[1:d]}\Psi(x) \\ \bar{\Psi}(\hat{K}_{[1:d]}\Psi(x))
	\end{pmatrix} = \hat{K}_{[1:d]}\Psi(x),
\end{align*}
where $\hat{K}_{[1:d]} \in \mathbb{R}^{d\times N}$ and $\hat{K}_{[d+1:N]} \in \mathbb{R}^{N-d\times N}$ are the first~$d$ and last~$N-d$ rows of~$\hat{K}$, respectively. Thus, only the first $d$ rows corresponding to the dynamics of the state are relevant for the predictions. Hence, Koopman-based prediction with coordinate projection resembles a discrete-time version of SINDy~\cite{BrunProc16} without the sparsity aspect.
% This shows that, interestingly, the last $N-d$ rows of $\hat{K}$ are completely removed by the coordinate projection.

\begin{figure}
	\centering
	\includegraphics[width=.7\linewidth]{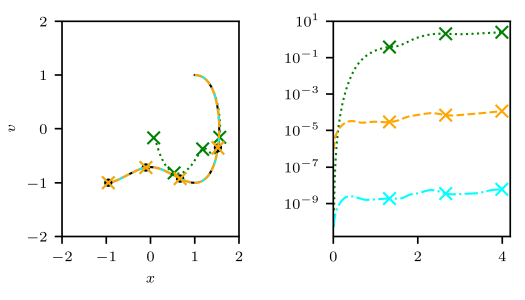}
	\caption{Comparison of the surrogate models with coordinate projection using $\mathbb{V}_3$ (\ref{eq:koopman_approx_dynamics}; \textcolor{orange}{$--$}), coordinate projection using $\mathbb{V}_5$ (\ref{eq:koopman_approx_dynamics}; \textcolor{cyan}{$\cdot-$}) and no projection using $\mathbb{V}_5$ (\ref{eq:projection_none}; \textcolor{green}{$\cdots$}) for Example~\ref{ex:duffing} on $\mathbb{X} = [-2,2]^2$ in comparison to the exact solution (---): exemplary trajectories (left) and mean error over time averaged over $x^0 \in \mathbb{X}$ (right). 
	}
	\label{fig:example_trajectory_duffing_order5}
\end{figure}

\noindent The benefits %importance 
of using such a projection~$\pi$ are demonstrated in the following example.
\begin{example}\label{ex:duffing}
	On the domain $\mathbb{X} = [-2, 2] \times [-2,2]$, we consider the unforced and undamped Duffing oscillator,
	\begin{align}\label{eq:duffing_example}
		\dot{x} = v,\qquad 
		\dot{v} = x - x^3.
	\end{align}
	For $n \in \mathbb{N}$, define the monomial dictionaries
	\begin{align}\label{eq:dictionary_monomial}
		\mathbb{V}_n = \operatorname{span} \{ x^a v^b \mid a,b \in [0:n]: a+b \leq n \}.
	\end{align}
	Using the time step~$\delta=0.01$, we approximate the Koopman operator as described in Section~\ref{sec:edmd}.
	The dynamics of the corresponding lifted, i.e., %standard 
	not projected, data-based surrogate model %in the Koopman framework approximated Koopman dynamics 
	are obtained simply by %simply integrating the% lifted system,
	\begin{align}\label{eq:projection_none}
		z(n+1) = \hat{K} z(n), \qquad z(0) = \Psi(x(0)).
	\end{align}
	Then, we have $x(n) = z_1(n), v(n) = z_2(n)$, $n \in \mathbb{N}_0$.
	Figure~\ref{fig:example_trajectory_duffing_order5} shows that the additional projection step in the dynamics~\eqref{eq:koopman_approx_dynamics} significantly improves the approximation accuracy and allows for predictions on much larger time intervals in comparison to its counterpart~\eqref{eq:projection_none} without projection step.
\end{example}
\noindent While the additional projection step typically yields a significant improvement of the approximation quality, we provide some additional insight why the coordinate projection is particularly well suited for the Duffing oscillator considered in Example~\ref{ex:duffing}. 
The reasoning is the close relationship between the system equations~\eqref{eq:duffing_example} and the choice of the dictionary~$\mathbb{V}_5$, which can be shown in a more general fashion: 
Let $\mathbb{V}$ be a dictionary including the coordinate functions, i.e., $\Psi$ can be written as~\eqref{eq:Psi}, and assume that $f_i \in \mathbb{V}$ holds for $i \in [1:d]$ (each component of the right hand side~\eqref{eq:ode} is contained in the dictionary).
The representation~\eqref{eq:Psi} allows to rewrite the argument of the regression problem~\eqref{eq:fitting} as
\begin{align*}
	& \|\Psi(x(t;\hat{x})) - K\Psi(\hat{x})\|^2_2 
	\!=\!  \left\| \begin{pmatrix} x(t;\hat{x}) -K_{[1:d]} \Psi(\hat{x}) \\ \bar{\Psi}(x(t;\hat{x})) - K_{[d+1:N]} \Psi(\hat{x}) \end{pmatrix} \right\|^2_2 
	\!\!=\!  \left\| \begin{pmatrix} \hat{x} + t f(\hat{x}) + \mathcal{O}(t^2) - K_{[1:d]} \Psi(\hat{x}) \\ \bar{\Psi}(x(t;\hat{x})) - K_{[d+1:N]} \Psi(\hat{x}) \end{pmatrix} \right\|^2_2\!\!,
\end{align*}
where we tacitly imposed sufficient smoothness of the vector field~$f$ such that $x(t;\hat{x}) = \hat{x} + t f(\hat{x}) + \mathcal{O}(t^2)$ holds.
Then, invoking $f_i \in \mathbb{V}$, the approximation error of the solution~$\hat{K}_{[1:d]}$ is bounded by $\mathcal{O}(t^2)$. % since $\hat{x}_i + tf_i(x) \in \mathbb{V}$ for each $i$. 
In particular, we obtain
\begin{align*}
	\hat{F}_{\pi}(\hat{x}) = x(t;\hat{x}) + \mathcal{O}(t^2)
\end{align*}
for all initial conditions $\hat{x} \in \mathbb{X}$.

\section{Closest-Point Projections}
\label{sec:closest_point_projection}
In the following we assume that the set~$M$ defined by~\eqref{eq:M}, which is induced by $\Psi \in L^2(\mathbb{X}, \mathbb{R})^N$, is a smooth and $d$-dimensional embedded manifold in $\mathbb{R}^N$, where continuous differentiability is a sufficient smoothness assumption for our purposes. Alternatively, one may consider any injective immersion $\Psi$ which satisfies any of the conditions of \cite[Proposition 4.22]{Lee12}.

Let $W \in \mathbb{R}^{N \times N}$ be a positive semi-definite matrix and define the weighted semi-inner product %sub-metric
\begin{align*}
	\langle u_1, u_2 \rangle_W := u_1^\top W u_2 \qquad\forall\,u_1,u_2 \in \mathbb{R}^N.
\end{align*}
If $W$ is invertible or the weaker condition
\begin{align}\label{eq:metric_determinant_condition}
	\det(\mathrm{D}\Psi(x)^\top W \mathrm{D}\Psi(x)) \neq 0 \qquad\forall\,x \in \mathbb{X}
\end{align}
holds, then $W$ induces a Riemannian metric on~$M$.
This then defines the notion of distance on the manifold~$M$. 

Based on the chosen semi-inner product~$\langle \cdot, \cdot \rangle_W$, we construct the \emph{closest-point projection}.
\begin{definition}\label{dfn:closest_point_projection}\ %
	For a given $(N \times N)$-matrix $W = W^\top \geq 0$ satisfying Condition~\eqref{eq:metric_determinant_condition}, the \textit{closest-point projection} $\pi_W : \mathbb{R}^N \to M$ is defined as
	\begin{align}\label{eq:projection_W}
		\pi_W(z) := \operatorname{argmin}_{p\in M} \| z - p\|_W.
	\end{align}
\end{definition}
\noindent Condition~\eqref{eq:metric_determinant_condition} ensures that $\pi_W$ is well-defined in a neighbourhood of the manifold embedded in $\mathbb{R}^N$, and, in particular, $\pi_W(z) = z$ for all $z \in {M}$.
Since the projection operator~\eqref{eq:projection_W} is invariant under scalings of $W$, i.e., $\pi_W(z) = \pi_{\alpha W}(z)$ holds for all $z \in \mathbb{R}^N$ and $\alpha > 0$, it suffices to consider semi-inner product~$W$ satisfying $\pdet(W)=1$. This corresponds to the choices of~$W$ for which the Riemannian volume of~$M$ is constant.

Next, we show that the coordinate projection is a closest-point projection. 
\begin{proposition}
	%The relationship with coordinate projection is less obvious, but can be seen as follows.
	Let $\Psi$ be given by Equation~\eqref{eq:Psi}, i.e., $\Psi$ contains the coordinate functions.
	Then, %the matrix $C \in \mathbb{R}^{N \times N}$ defined by
	\begin{align}
		C := \begin{pmatrix}
			I_d & 0_{d\times N-d} \\
			0_{N-d \times d} & 0_{N-d \times N-d}
		\end{pmatrix} \in \mathbb{R}^{N \times N}
	\end{align}
	induces a Riemannian metric on ${M}$, and the coordinate projection~\eqref{eq:coordinate_projection} coincides with the closest-point projection with $W = C$.  
\end{proposition}
\begin{proof}
	We verify Condition~\eqref{eq:metric_determinant_condition} %in order 
	to show that $C$ induces a Riemannian metric on ${M}$. Let $x \in \mathbb{X}$ be given. Then, we have $\det(\mathrm{D}\Psi(x)^\top C \mathrm{D}\Psi(x)) = \det(I_d) \neq 0$, which can be directly inferred by rewriting the left hand side as
	\begin{align*}
		&\det \left( \begin{pmatrix}
			I_d  \\ \mathrm{D} \bar{\Psi}(x)
		\end{pmatrix}^\top \begin{pmatrix}
			I_d & 0_{d\times N-d} \\
			0_{N-d \times d} & 0_{N-d \times N-d}
		\end{pmatrix} \begin{pmatrix}
			I_d  \\ \mathrm{D} \bar{\Psi}(x)
		\end{pmatrix} \right).
	\end{align*}
	Now, for any $z \in \mathbb{R}^N$, one has
	\begin{align*}
		\pi_C(z) &= \operatorname{argmin}_{p\in M} \| z - p\|_C \\
		&= \Psi \circ \operatorname{argmin}_{x \in \mathbb{X}} \| z - \Psi(x)\|_C \\
		&= \Psi \circ \operatorname{argmin}_{x \in \mathbb{X}} \left\|  \begin{pmatrix}
			z_{[1:d]} -  x \\ z_{[d+1:N]} - \bar{\Psi}(x)
		\end{pmatrix} \right\|_C \\
		&= \Psi \circ \operatorname{argmin}_{x \in \mathbb{X}} \left\| z_{[1:d]} - x \right\| = \Psi (z_{[1:d]}),
	\end{align*}
	which completes the proof.
\end{proof}
We emphasize that the matrix~$C$ is not invertible. Further, the closest-point projection~$\pi_W$ is, in general, a nonlinear projection, and, may be implemented using variants of steepest descent or Newton's method~\cite{AbsiMaho09}.

The following proposition shows that the error resulting from the closest-point projection is proportionally bounded to the approximation error in the regression problem~\eqref{eq:fitting}.
\begin{proposition}\label{prop:projbound}
	The closest point projection of Definition~\ref{dfn:closest_point_projection} satisfies
	\begin{align*}
		\left\Vert \pi_W(\hat{K}\Psi(x)) - \mathcal{K}^t\Psi(x) \right\Vert_W
		\leq 2 \left\Vert \hat{K} \Psi(x) - \mathcal{K}^t \Psi(x) \right\Vert_W
	\end{align*}
	for all $x\in \mathbb{R}^d$.
	That is, its error is bounded by twice the training error in the given metric $W$.
\end{proposition}
\begin{proof}
	Since $\mathcal{K}^t \Psi(x) \in M$, the definition of $\pi_W(\hat{K}\Psi(x))$ implies
	\begin{align*}
		\left\Vert \pi_W(\hat{K} \Psi(x)) - \hat{K} \Psi(x) \right\Vert_W
		\leq \left\Vert \mathcal{K}^t \Psi(x) - \hat{K} \Psi(x) \right\Vert_W.
	\end{align*}
	Hence, using the triangle inequality and the definition of the regression problem~\eqref{eq:fitting} weighted with $W$ yields the assertion. 
\end{proof}

One may observe that the bound of Proposition~\ref{prop:projbound} uses the semi-inner product~$W$, which is not the metric used in the construction of the surrogate model~$\hat{K}$, cp.\ the regression problem~\eqref{eq:fitting}. However, we show in the following proposition that the solution of problem~\eqref{eq:fitting} also solves the weighted regression problem.
\begin{proposition}
	The solution $\hat{K}$ to the regression problem~\eqref{eq:fitting} satisfies %also solves the weighted regression problem, i.e.,
	\begin{align}\label{eq:fitting_weighted}
		\hat{K}
		& \in \operatorname{arg} \hspace*{-0.5cm}\min_{K\in \mathbb{R}^{N \times N} \hspace*{0.4cm}}\hspace*{-0.2cm}
		\int_\mathbb{X} \|\Psi(x(t;\hat{x})) - K\Psi(\hat{x})\|^2_W \,\mathrm{d}\hat{x}
	\end{align}
	for all $W = W^\top \geq 0$, $W \in \mathbb{R}^{N \times N}$.
\end{proposition}
\begin{proof}
	Using $\Psi_X,\Psi_Y$ given by~\eqref{eq:KI}, differentiation of the cost w.r.t.\ ${K}$ in an arbitrary direction $\Delta$ yields
	\begin{align*}
		2 \int_\mathbb{X} \langle K\Psi(\hat{x}) - \Psi(x(t;\hat{x})) \,,\, \Delta\Psi(\hat{x}) \rangle_W \,\mathrm{d}\hat{x} 
		= & 2 \langle \int_\mathbb{X} K \Psi(\hat{x}) \Psi(\hat{x})^\top - \Psi(x(t; \hat{x})) \Psi(\hat{x})^\top \,\mathrm{d}\hat{x} , \Delta \rangle_{W} \\
		= & 2 \langle K \Psi_X - \Psi_Y, \Delta \rangle_W = 2 \langle W (K \Psi_X - \Psi_Y), \Delta \rangle.
	\end{align*}
	This expression is nullified for $\hat{K} = \Psi_Y \Psi_X^{-1}$, which is thus indeed a solution to the regression problem~\eqref{eq:fitting_weighted}; independently of the chosen weighting matrix $W = W^\top \geq 0$. %This completes the proof.
\end{proof}

\section{Geometric Projection}\label{sec:geometric_projection}
The choice of a suitable metric is a non-trivial task. 
We propose a geometrically-motivated one, which exhibits a superior performance as demonstrated in this section. % example.
\begin{definition}
	Let $\Sigma \in \mathbb{R}^{N\times N}$ be given by
	\begin{align} \label{eq:sigma_definition}
		\Sigma = \int_{\mathbb{X}} (\hat{K} \Psi(x) - \mathcal{K}^t \Psi(x))(\hat{K} \Psi(x) - \mathcal{K}^t \Psi(x))^\top \,\mathrm{d}x
	\end{align}
	and be invertible. Then, the \emph{geometric projection}~$\pi_W$ is the closest point projection of Definition~\ref{dfn:closest_point_projection} associated with the metric $W = \det(\Sigma)^{1/N} \Sigma^{-1}$. %\in \operatorname{met}^1(M)$.
\end{definition}
Recalling the scaling invariance, it is straightforward to see that the geometric projection is a special case of closest-point projection with the metric $W = \Sigma^{-1}$, where $\Sigma$ is defined as in \eqref{eq:sigma_definition}.
The geometric projection can be interpreted from a probabilistic viewpoint. 
Suppose that we approximate the Koopman operator based on normally distributed i.i.d.\ random variables, i.e., $\hat{K} \Psi(x) \sim N(\mathcal{K}^t \Psi(x), \Sigma)$ for every $x \in D$ with~$\Sigma$ defined by~\eqref{eq:sigma_definition}.
Then, for every $x \in D$, the likelihood of a point $p \in \mathbb{R}^N$ being equal to $\mathcal{K}^t \Psi(x)$ can be written as
\begin{align}\nonumber
	\rho(\,p\, |\, \hat{K} \Psi(x)) = \frac{\exp\left( -0.5 \cdot \| p - \hat{K} \Psi(x) \|_{\Sigma^{-1}}^2 \right)}{\sqrt{(2 \pi)^N \det(\Sigma)}}.
\end{align}
If we restrict ourselves to look for points $p \in M$, then the maximum likelihood solution is exactly
\begin{align*}
	\operatorname{argmax}_{p\in M} \rho(\,p |\, \hat{K} \Psi(x))
	= \pi_W(\hat{K} \Psi(x)),
\end{align*}
where $W = \Sigma^{-1}$ holds.

Next, we consider the example of a pendulum to %numerically 
demonstrate the advantages of the geometric projection in comparison to its coordinate-based counterpart. %the coordinate projection. 
To this end, we use the notation~$\Delta t$ in the following to indicate the time step, i.e., $\hat{K}$ approximates $\mathcal{K}^{\Delta t}$.
\begin{example}\label{ex:pendulum}
	Consider the pendulum with dynamics %given by
	\begin{align}\label{eq:pendulum_example}
		\dot{x} = v, \qquad \dot{v} = -\sin(x)
	\end{align}
	on the domain $\mathbb{X} = [-\pi, \pi] \times [-3,3]$. We approximate the Koopman operator by taking $10,000$ data points~$\hat{x}$ drawn uniformly i.i.d.\ from~$\mathbb{X}$ and the respective solution $x(\Delta t; \hat{x})$ as described in Section~\ref{sec:edmd} using monomial dictionaries~$\mathbb{V}_n$, $n \in \mathbb{N}$, cp.~\eqref{eq:dictionary_monomial}.
\end{example}

\begin{figure}[!htb]
	\centering
	\includegraphics[width=.7\columnwidth]{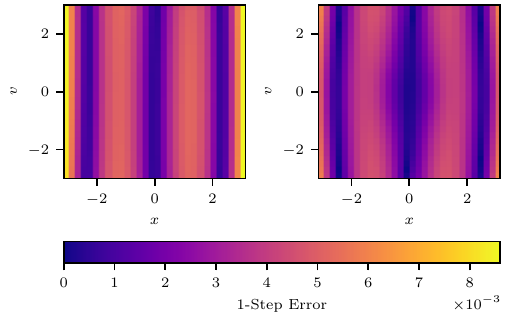}
	\caption{Example~\ref{ex:pendulum} with dictionary~$\mathbb{V}_2$: Comparison of the one-step errors~\eqref{eq:1step} for the approximated dynamics using coordinate (left) and geometric projection (right), respectively.
	}
	\label{fig:pendulum_1step_heat_order3}
\end{figure}

\begin{figure}[!htb]
	\centering
	\includegraphics[width=.45\columnwidth]{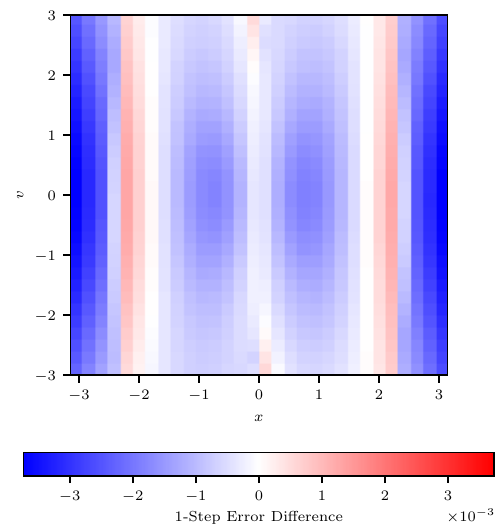}
	\includegraphics[width=.45\columnwidth]{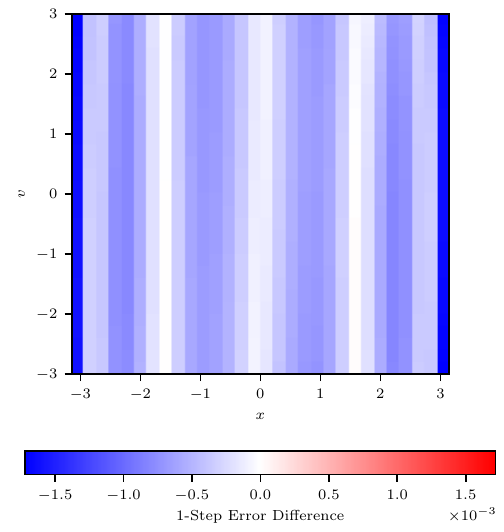}
	\caption{The difference in 1-Step errors $E_{1-\text{Step}}^W(\hat{x}) - E_{1-\text{Step}}^C(\hat{x})$ for the approximated dynamics using geometric and coordinate projections and dictionary orders 2 (left) and 3 (right).}
	\label{fig:pendulum_1step_heat}
\end{figure}

The projection methods are compared by examining the approximated system dynamics~\eqref{eq:koopman_approx_dynamics} with those obtained through a high-order numerical integration scheme with step-size control. The \textit{one-step error} at a given value $\hat{x} \in \mathbb{X}$ is %computed as
\begin{align}\label{eq:1step}
	E_{1-\text{Step}}^W(\hat{x}) := \| \hat{F}_{\pi_W}(\hat{x}) - x(\Delta t;\hat{x}) \|
\end{align}
for a given projection~$\pi_W$.
%, where $x(\delta;x^0)$ is obtained using a fourth-order Runge-Kutta method.
Figure~\ref{fig:pendulum_1step_heat_order3} shows the one-step errors for the coordinate and geometric projections using the dictionary~$\mathbb{V}_2$.
%\vspace*{-1.5cm}
For a better comparison and the impact of using, in addition, monomials of order three in the dictionary, Figure~\ref{fig:pendulum_1step_heat} shows the difference in one-step errors. 
%\vspace*{-0.5cm}
For~$\mathbb{V}_2$, the geometric projection has a lower one-step error in most of the domain, which is consistent with the underlying regression problem, in which the $L^2$-error is minimized. For~$\mathbb{V}_3$, the geometric projection has a lower one-step error everywhere. %As the dictionary size is increased, the projection error decreases~\cite{SchaWort22}. Here, 
The underlying reasoning is that the geometric projection exploits its additional freedom by taking more dictionary elements into account to return to the manifold~$M$.

\begin{figure}[!htb]
	\centering
	\includegraphics[width=0.49\linewidth]{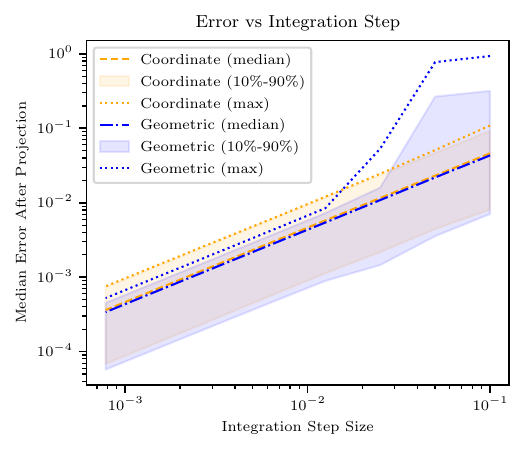}
	\includegraphics[width=0.49\linewidth]{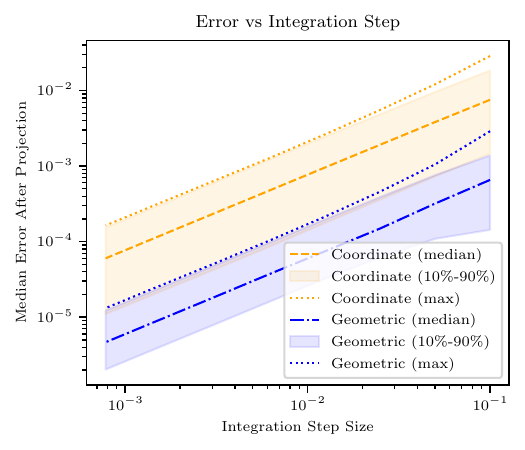}
	\caption{Comparison of median 1-step error for the geometric and coordinate projections using~$\mathbb{V}_2$ (left) and~$\mathbb{V}_3$ (right) depending on time step~$\Delta t$.
	}
	\label{fig:timestep_errors}
\end{figure}

The 1-step error depends on $\Delta t$ for any projection method. 
Figure~\ref{fig:timestep_errors} shows how the 1-step error statistics change for the coordinate and geometric projections depending on the time step~$\Delta t$.
Again, the geometric projection clearly outperforms its coordinate-based counterpart if the dictionary size increases, i.e., for $\mathbb{V}_n$ and $n \in \{3,4,5\}$ (see Figure~\ref{fig:timestep_errors} for $n=2,3$). %Both projections are improved for all times by the addition of third order dictionary terms.
%However, for the third-order dictionary $\mathbb{V}_3$ the geometric projection also clearly outperforms the coordinate projection.
Regardless of dictionary size, both projections degrade as the time-step is increased.

\noindent Next, we consider the Lorenz system for a dictionary~$\Psi$ without $\psi_1(x,y,z) = x$, i.e., the dictionary does not contain all coordinate functions, see also \cite[Section III.C]{li2017extended}.
\begin{example}\label{ex:lorenz}%[Lorenz Attractor]
	Consider the Lorenz system given by %attractor dynamics,
	\begin{align}\nonumber %\label{eq:lorenz_attractor}
		\dot{x} &= \sigma (y-x), &
		\dot{y} &= x(\rho-z)-y, &
		\dot{z} &= xy - \beta z,
	\end{align}
	on the domain $\mathbb{X} = [-20,20]\times[-20,20]\times[10,50]$.
\end{example}

\begin{figure}[htb]
	\centering
	\includegraphics[width=.7\linewidth,trim={0 0.1cm 0 0.3cm},clip]{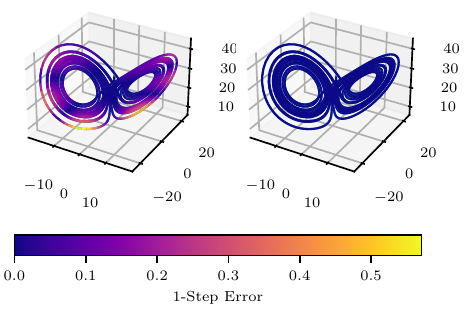}
	\caption{1-step errors along a sample trajectory of %the Lorenz system 
		Example~\ref{ex:lorenz} for the coordinate (left) and geometric projections (right) using the dictionary $\mathbb{V}_4 \setminus \{ x \}$.}
	\label{fig:lorenz_example_trajectories}
\end{figure}

%We examine the effect of excluding the coordinate function  from the dictionary.
\noindent In this case, the coordinate projection can be realised by selecting three observables from which $x,y,z$ can be %theoretically 
recovered: $xz$, $y$, and $z$, i.e., $x$ can be reconstructed by $x = xz/z$ since $z \in [10, 50]$.
Figure \ref{fig:lorenz_example_trajectories} shows an examplary trajectory from the true and approximated dynamics using the dictionary %ies $\mathbb{V}_4$ and 
$\mathbb{V}_4 \setminus \{ x \}$.
While the geometric projection performs well without the coordinate function, the accuracy of the reconstructed coordinate projection is poor.

\section{Conclusions and outlook}\label{sec:conclusions}
Our key contributions are the following: 
First, we have demonstrated the need for a reprojection step whenever the dictionary is not Koopman invariant. 
Second, we proposed a general framework to conduct the reprojection step based on a large class of semi-inner products. In particular, it is not necessary that the coordinate functions are contained in the dictionary to conduct the reprojection step, see, e.g., the novel closest-point projection supposing invertibility of the map~$\Psi$. 
Third, we have rigorously shown that the additional reprojection step essentially maintains the estimation error resulting from the regression problem and, indeed, significantly reduces the approximation error as shown in our numerical simulations. A key reason is that the chosen weighting does not interfere with the regression problem to be solved for computing the data-based compression~$\hat{K}$.

Clearly, the proposed framework is directly applicable to nonlinear control-affine systems, if bilinear surrogate models are used, see, e.g., \cite{otto2021koopman,SchaWort22} and the references therein.
Here, already the coordinate projection has turned out to be very beneficial in simulation and experiments for non-holonomic robots~\cite{BoldEsch23} such that we expect clear benefits if the novel geometric projection is applied.
This claim also applies to eDMD with control~\cite{proctor2018generalizing} since the weighting~$W$ of the control in the augmented state may be set to zero.

Future work might be devoted to leveraging recently introduced concepts~\cite{Mezi20} for the analysis of systems with a globally-stable attractor in our setting more tailored towards control systems typically lacking such structures. Here, taking into account the recent results of~\cite{GoswPale21} is of interest. Furthermore, we will leverage known results from regression~\cite{Bish06} to further analyze and potentially improve the proposed framework. 

\bibliographystyle{abbrv}
\bibliography{ref}

\begin{thebibliography}{10}

\bibitem{AbsiMaho09}
P.-A. Absil, R.~Mahony, and R.~Sepulchre.
\newblock Optimization algorithms on matrix manifolds.
\newblock In {\em Optimization Algorithms on Matrix Manifolds}. Princeton
  University Press, 2009.

\bibitem{BevaSosnHirc21}
P.~Bevanda, S.~Sosnowski, and S.~Hirche.
\newblock {K}oopman operator dynamical models: Learning, analysis and control.
\newblock {\em Annu Rev Control}, 52:197--212, 2021.

\bibitem{Bish06}
C.~M. Bishop and N.~M. Nasrabadi.
\newblock {\em Pattern recognition and machine learning}, volume~4.
\newblock Springer, 2006.

\bibitem{BoldEsch23}
L.~Bold, H.~Eschmann, M.~Rosenfelder, H.~Ebel, and K.~Worthmann.
\newblock On {K}oopman-based surrogate models for non-holonomic robots.
\newblock {\em Preprint arxiv:2303.09144}, 2023.

\bibitem{BoldGrun23}
L.~Bold, L.~Grüne, M.~Schaller, and K.~Worthmann.
\newblock Practical asymptotic stability of data-driven model predictive
  control using extended {DMD}, 2023.
\newblock Preprint available on arXiv.

\bibitem{BrunBrun16}
S.~L. Brunton, B.~W. Brunton, J.~L. Proctor, and J.~N. Kutz.
\newblock {K}oopman invariant subspaces and finite linear representations of
  nonlinear dynamical systems for control.
\newblock {\em PloS one}, 11(2):e0150171, 2016.

\bibitem{BrunProc16}
S.~L. Brunton, J.~L. Proctor, and J.~N. Kutz.
\newblock Discovering governing equations from data by sparse identification of
  nonlinear dynamical systems.
\newblock {\em Proceedings of the national academy of sciences},
  113(15):3932--3937, 2016.

\bibitem{GoswPale21}
D.~Goswami and D.~A. Paley.
\newblock Bilinearization, reachability, and optimal control of control-affine
  nonlinear systems: A {K}oopman spectral approach.
\newblock {\em IEEE Transactions on Automatic Control}, 67(6):2715--2728, 2021.

\bibitem{KlusNusk20}
S.~Klus, F.~N{\"u}ske, S.~Peitz, J.-H. Niemann, C.~Clementi, and
  C.~Sch{\"u}tte.
\newblock Data-driven approximation of the {K}oopman generator: Model
  reduction, system identification, and control.
\newblock {\em Physica D: Nonlinear Phenomena}, 406:132416, 2020.

\bibitem{KordMezi18b}
M.~Korda and I.~Mezi{\'c}.
\newblock Linear predictors for nonlinear dynamical systems: {K}oopman operator
  meets model predictive control.
\newblock {\em Automatica}, 93:149--160, 2018.

\bibitem{KordMezi18}
M.~Korda and I.~Mezi{\'c}.
\newblock On convergence of extended dynamic mode decomposition to the
  {K}oopman operator.
\newblock {\em J. Nonlinear Sci.}, 28(2):687--710, 2018.

\bibitem{Lee12}
J.~M. Lee.
\newblock {\em Introduction to Smooth manifolds}.
\newblock Springer, 2012.

\bibitem{li2017extended}
Q.~Li, F.~Dietrich, E.~M. Bollt, and I.~G. Kevrekidis.
\newblock Extended dynamic mode decomposition with dictionary learning: A
  data-driven adaptive spectral decomposition of the {K}oopman operator.
\newblock {\em Chaos: An Interdisciplinary J. Nonlinear Sci.}, 27(10):103111,
  2017.

\bibitem{MaurGonc16}
A.~Mauroy and J.~Goncalves.
\newblock Linear identification of nonlinear systems: A lifting technique based
  on the {K}oopman operator.
\newblock In {\em 55th IEEE Conference on Decision and Control (CDC)}, pages
  6500--6505, 2016.

\bibitem{MaurMezi20}
A.~Mauroy, I.~Mezi{\'c}, and Y.~Susuki.
\newblock {\em The {K}oopman operator in systems and control}.
\newblock Springer Nature, Cham, Switzerland, Feb. 2020.

\bibitem{Mezi20}
I.~Mezi{\'c}.
\newblock Spectrum of the {K}oopman operator, spectral expansions in functional
  spaces, and state-space geometry.
\newblock {\em Journal of Nonlinear Science}, 30(5):2091--2145, 2020.

\bibitem{NuskPeit23}
F.~N{\"u}ske, S.~Peitz, F.~Philipp, M.~Schaller, and K.~Worthmann.
\newblock Finite-data error bounds for {K}oopman-based prediction and control.
\newblock {\em J. Nonlinear Sci.}, 33:14, 2023.

\bibitem{otto2021koopman}
S.~E. Otto and C.~W. Rowley.
\newblock {K}oopman operators for estimation and control of dynamical systems.
\newblock {\em Annu. Rev. Control Robot. Auton. Syst.}, 4:59--87, 2021.

\bibitem{PhilScha23}
F.~Philipp, M.~Schaller, K.~Worthmann, S.~Peitz, and F.~N{\"u}ske.
\newblock Error bounds for kernel-based approximations of the {K}oopman
  operator, 2023.
\newblock Submitted, Preprint: arXiv:2301.08637.

\bibitem{proctor2018generalizing}
J.~L. Proctor, S.~L. Brunton, and J.~N. Kutz.
\newblock Generalizing {K}oopman theory to allow for inputs and control.
\newblock {\em SIAM J. Appl. Dyn. Syst.}, 17(1):909--930, 2018.

\bibitem{SchaWort22}
M.~Schaller, K.~Worthmann, F.~Philipp, S.~Peitz, and F.~N{\"u}ske.
\newblock Towards reliable data-based optimal and predictive control using
  extended {DMD}.
\newblock {\em IFAC-PapersOnLine}, 56(1):169--174, 2023.

\bibitem{TeraShir21}
H.~Terao, S.~Shirasaka, and H.~Suzuki.
\newblock Extended dynamic mode decomposition with dictionary learning using
  neural ordinary differential equations.
\newblock {\em Nonlinear theory appl.\ IEICE}, 12(4):626--638, 2021.

\bibitem{WillHema2016}
M.~O. Williams, M.~S. Hemati, S.~T. Dawson, I.~G. Kevrekidis, and C.~W. Rowley.
\newblock Extending data-driven {K}oopman analysis to actuated systems.
\newblock {\em IFAC-PapersOnLine}, 49(18):704--709, 2016.

\bibitem{WillMatt15}
M.~O. Williams, I.~G. Kevrekidis, and C.~W. Rowley.
\newblock A data--driven approximation of the {K}oopman operator: Extending
  dynamic mode decomposition.
\newblock {\em J. Nonlinear Sci.}, 25:1307--1346, 2015.

\bibitem{ZhanZuaz23}
C.~Zhang and E.~Zuazua.
\newblock A quantitative analysis of {K}oopman operator methods for system
  identification and predictions.
\newblock {\em Comptes Rendus. M{\'e}canique}, 351(S1):1--31, 2023.

\end{thebibliography}

\end{document}